\documentclass[12pt]{article}
\usepackage{amsmath, amsthm, amssymb, amscd, graphicx}

\usepackage{geometry} 

\title{Interval exchanges that do not embed in free groups}
\author{Christopher F. Novak}
\date{July 22, 2010} 

\newtheorem{theorem}{Theorem}[section]

\newtheorem{lemma}[theorem]{Lemma}
\newtheorem{proposition}[theorem]{Proposition}

\theoremstyle{remark}
\newtheorem{remark}[theorem]{Remark}

\theoremstyle{definition}

\theoremstyle{question}
\newtheorem{question}[theorem]{Question}

\def\co{\colon\thinspace}

\begin{document}

\maketitle

\begin{abstract}
A disjoint rotation map is an interval exchange transformation (IET) on the unit interval that
acts by rotation on a finite number of invariant subintervals. It is currently unknown 
whether the group $\mathcal{E}$ of all IETs possesses any non-abelian free subgroups. It is shown
 that it is not possible for a disjoint rotation map to occur in a subgroup of $\mathcal{E}$ that is isomorphic
 to a non-abelian free group.   
 \end{abstract}

\section{Introduction}

An interval exchange transformation (IET) is an invertible map $[0,1) \rightarrow [0,1)$   
defined by a finite partition of $[0,1)$ into half-open subintervals and 
 a reordering of these intervals by translation. The dynamics of single IETs have been actively studied 
 since the late 1970s. See the recent survey of Viana {\bf \cite{Viana06}} for a presentation of many early results in this area.
 The study of IET dynamics is currently quite active, due in part to the close connection between IETs and the moduli space of translation surfaces (see 
 the survey of Zorich {\bf\cite{Zorich06}}), and also due to the recent resolution of some long-outstanding problems in the area (e.g.,  {\bf\cite{AvilaForni07}}, 
 {\bf\cite{AvilaViana07}}, and {\bf\cite{FerencziHoltonZamboni04}}).
      
More recently, group actions by interval exchanges have begun to be studied; see, for instance, {\bf \cite{Arnoux81}}, {\bf \cite{Novak09}}, and {\bf \cite{Novak10}}. 
The set $\mathcal{E}$ of all interval exchanges forms a group under composition, 
and an \emph{interval exchange action} of a group $G$ is a homomorphism $G \rightarrow \mathcal{E}$. The study of such actions is motivated by 
the analogous study of group actions on manifolds, particularly 1--dimensional ones, by means of homeomorphisms or diffeomorphisms. The subject of group 
actions on 1--manifolds is well-developed and quite active; see, for instance, {\bf \cite{Ghys01}} or {\bf \cite{Navas07bookEngArxiv}}--{\bf \cite{Navas07book}}. 

In contrast, many fundamental questions which are well understood for groups acting on 1--manifolds 
are currently open for the group of IETs. Perhaps foremost among
these is the following question, due to A. Katok:

\begin{question} \label{FreeIETgroupsQ}
 Does $\mathcal{E}$ contain a subgroup isomorphic to $F_2$, the non-abelian free group on two generators?
\end{question}

It is easy to construct examples of non-abelian free subgroups in $\text{Diff}(S^1)$ and $\text{Diff}(\mathbb{R})$ by means of the ping-pong construction. 
 More detailed results, analogous to the Tits' alternative, are known for $\text{Homeo}_+(S^1)$ and $\text{Diff}_+^{\, \omega}(S^1);$ see
  {\bf \cite{Margulis00}}  and {\bf \cite{FarbShalen02}}, respectively. It is also shown in {\bf \cite{Ghys01}} that for a residual set of pairs $(f,g)$ in
   $\text{Homeo}_+(S^1)$, the group $\langle f,g \rangle$ is isomorphic to $F_2$. 
   However, there are examples of groups of homeomorphisms of 1--manifolds that do not contain non-abelian
  free subgroups. For instance, it is known from work of Brin and Squier {\bf \cite{BrinSquier85}} that this is the case for the group $\text{PL}_+([0,1])$ of piecewise-linear 
  homeomorphisms of the interval. 
  
 \begin{remark}
   The paper {\bf \cite{BrinSquier85}} also shows that the mechanism by which $\text{PL}_+([0,1])$ fails to contain non-abelian free subgroups is not
   an obvious one. In particular, this work proves that $\text{PL}_+([0,1])$ does not satisfy a law; i.e., there 
   does \emph{not} exist $\omega \in F_2 \setminus \{e\}$ such that $\phi(\omega) = id$ for every homomorphism $\phi \co \! F_2 \rightarrow \text{PL}_+([0,1]).$  
   It is also the case that $\mathcal{E}$ does not satisfy a law. If such a law were to exist, then it would have to be satisfied
    by every finite group, since every finite group occurs as a subgroup of $\mathcal{E}$. The existence of such a universal law for finite groups 
    is impossible; in particular, it would imply that $F_2$ is not residually finite, which is false (see section III.18 of {\bf \cite{delaHarpe00}}). 
 \end{remark} 
  
The current work shows that a particular class of IETs, the \emph{disjoint rotation maps}, cannot occur in a non-abelian
 free subgroup of $\mathcal{E}$, if such subgroups actually do exist.    
 Briefly, a disjoint rotation map $r$ is an IET for which there is a finite partition of $[0,1)$ into  $r$-invariant subintervals $I_j,$ such that
 $r$ restricted to each $I_j$ is an exchange of two further subintervals; a graphical depiction is given in 
Figure \ref{Fig disjoint rotation map}, and a formal definition is provided below.

\begin{theorem} \label{disjoint rotation relation theorem}
Let $r$ be conjugate in $\mathcal{E}$ to a disjoint rotation map, and let $g\in \mathcal{E}$ be an arbitrary interval exchange. 
Then the subgroup $\langle r, g \rangle$ is not isomorphic to the free group on two generators.    
\end{theorem} 

The proof of Theorem \ref{disjoint rotation relation theorem} is essentially constructive, in that it describes a nontrivial word in the generators 
$r$ and $g$ that forms the identity map. The construction strongly relies on two features of the disjoint rotation map $r$; such maps have iterates that are arbitrarily 
close to the identity in an $L^1$ sense, and the iterates $r^n$ have essentially the same number of discontinuities as $r$ itself. It is known 
({\bf \cite{Veech84a}}, Theorems 1.3 and 1.4) that almost every irreducible
IET possesses iterates that are $L^1$ close to the identity, which raises the question of whether the construction described below can be adapted to such maps. The immediate 
obstruction is that for most such IETs $f$, the number of discontinuities of $f^n$ grows linearly with $n$; this prevents one from showing that 
a particular word formed from $f^n$ and $g$ has support contained in a neighborhood of a fixed finite set for infinitely many $n$.
   
It is interesting to note that conjugates of disjoint rotation maps are precisely those IETs that can occur 
in the image of a continuous homomorphism $\mathbb{R} \rightarrow \mathcal{E};$ see {\bf\cite{Novak10}}. 
Thus, Theorem \ref{disjoint rotation relation theorem} and the Tits' Alternative {\bf \cite{Tits72}}
imply that any linear Lie group that continuously embeds in $\mathcal{E}$ must be virtually solvable.

\subsubsection*{Acknowledgements}
The author is indebted to Professor Michael Lachance for his encouragement in this work.

\section{Proof of Theorem \ref{disjoint rotation relation theorem}}

\subsection{Notation}	
	We now give a precise definition and notation for interval exchanges.
Let $\pi \in \Sigma_n$ be a permutation of $\{1,2, \ldots, n\}$, and let $\lambda$ be a vector in the simplex 
\[ \Lambda_n = \left\lbrace \lambda = (\lambda_1, \ldots, \lambda_n)\co \lambda_i > 0, \
  \sum \lambda_i = 1     \right\rbrace \subseteq \mathbb{R}^n.  \]
   The vector $\lambda$ induces a partition of $[0,1)$ into intervals
  \begin{equation} \label{IET_form1} I_j = \left[ \beta_{j-1} := \sum_{i=1}^{i = j-1} \lambda_i,\ \, \beta_j := 
  \sum_{ i=1}^{i = j} \lambda_i \right), \ \ 1 \leq j \leq n. 
  \end{equation}  
Let $f_{(\pi, \lambda)}$ be the IET that translates each $I_j$ so that the ordering of these intervals
within $[0,1)$ is permuted according to $\pi$. More precisely, 
\begin{equation} \label{IET_form3} f_{(\pi, \lambda)} (x) = x + \omega_j, \quad \text{if }x \in I_j, \end{equation} 
where 
\begin{equation} \label{IET_form2} \omega_j = \Omega_\pi(\lambda)_j\, = \sum_{i:\, \pi(i) < 
\pi(j)} \lambda_i\ -\ \sum_{i:\, i<j} \lambda_i.  \end{equation}
Note that $\Omega_\pi \co \Lambda_n \rightarrow \mathbb{R}^n$ is a linear map depending only on $\pi$.

The above notation is adapted in the following way to represent a disjoint rotation map.
 Given $\lambda \in \Lambda_n$ for 
some $n$, define the points $\beta_j$  and the intervals $I_j$ by equation $\eqref{IET_form1}.$ 
Let $\alpha \in \mathbb{T}^n = (\mathbb{R} / \mathbb{Z})^n$ be given, where $\mathbb{T}^n$ is to be identified with $[0,1)^n$.

\begin{figure}[htb]
\begin{center}
\includegraphics{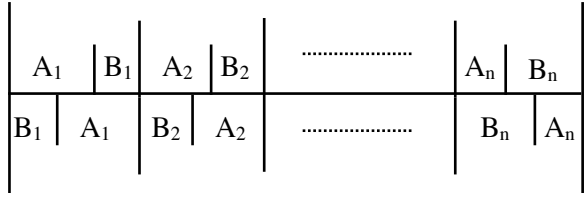}
\caption{A disjoint rotation map   \label{Fig disjoint rotation map}  }
\end{center}
\end{figure}

 Define the 
\emph{disjoint rotation map} $r_{[\alpha, \lambda]}$ by 

\begin{equation} \label{disjoint rotation map defn} r_{[\alpha, \lambda]}(x) = \left\lbrace \begin{array}{ll}
 x+\lambda_j {\alpha}_j, & x\in [\beta_{j-1},\, \beta_j - \lambda_j{\alpha}_j ) \\          
 x+ \lambda_j{\alpha}_j - \lambda_j, 
 & x\in [\beta_j -\lambda_j{\alpha}_j,\, \beta_j ). \end{array} \right.    \end{equation}	
See Figure \ref{Fig disjoint rotation map} for a graphical representation of a disjoint rotation map.
	
\subsection{Construction of a relation in $\langle r, g \rangle$}	
		
To begin the proof of Theorem \ref{disjoint rotation relation theorem}, it can be assumed after a conjugacy in $\mathcal{E}$ that the map $r$ has the 
form in equation \eqref{disjoint rotation map defn} for some $\alpha \in \mathbb{T}^n$ and $\lambda \in \Lambda_n$. If all points are periodic under $r$, then 
$r$ has finite order, and Theorem \ref{disjoint rotation relation theorem} holds trivially. Thus, assume that $r$ has infinite order; equivalently, assume 
$\alpha \notin (\mathbb{Q} / \mathbb{Z})^n$. If the set of periodic points  
of $r$ is nonempty, then after replacing $r$ by an 
iterate it can be assumed that all periodic points of $r$ are fixed. After a further conjugacy and possibly redefining $n$, it can be assumed that 
$\text{Fix}(r) = I_n$; in this case, $\alpha = (\alpha_1, \ldots, \alpha_{n-1}, 0)$ where $\alpha_i \in [0,1)$ is irrational for $1 \leq i \leq n-1$. Define 
the $\emph{support}$ of $r$, denoted $\text{supp}(r)$, to be the complement of its fixed point set.

A relation in $\langle r, g \rangle$ is constructed using the map $h = [g^{-1}, s^{-1}]\circ[g^{-1}, s],$  where $s = r^M$ for some integer $M$ to 
be chosen later and $[x, y] = xyx^{-1}y^{-1}$. It is to be shown that for suitably chosen $s = r^M$, the support of $h$ is
contained in a small neighborhood of a finite set. 

In particular, define the finite set $P$ by 
\[ P = \lbrace \beta_i \co 0 \leq i <  n   \rbrace \cup \lbrace g^{-1}(\beta_i) \co 0 \leq i < n  \rbrace \cup 
	\lbrace x \co g\text{ is discontinuous at $x$}\rbrace.   \]
Let $P' = P \cap\, \text{supp}(r)$. Since by assumption all non-fixed orbits of $r$ are infinite, it is possible to choose
 an integer $d > 0$ such that $r^d(P') \cap P' = \emptyset$. 
For $\epsilon >0$ and $p \in [0,1) \cong \mathbb{R} / \mathbb{Z},$ let $N_\epsilon(p)$ denote the open $\epsilon$-ball 
centered at $p$ in $\mathbb{R}/ \mathbb{Z}$; define the sets $X = X_\epsilon =  \bigcup_{p \in P} N_\epsilon(p)$ and $X' = (X \cap \text{supp}(r)).$
Next, fix $\epsilon > 0$ sufficiently small so that:

\begin{itemize}
	\item[i)] the collection of sets $\lbrace N_\epsilon(p) \co p \in P \rbrace$ are pairwise disjoint, and 
		
	\item[ii)]  the sets $X'$ and $r^d(X')$ are disjoint.

\end{itemize}

Finally, choose an integer $M>0$ such that the rotation rates $M\alpha_i \in \mathbb{R}/ \mathbb{Z}$ of $s = r^M$ are
 in $N_{\epsilon / 10}(0)$ for $i = 1, \ldots, n$. To see that 
such an $M$ exists, associate to $r$ the translation $\widehat{r}\co \mathbb{T}^n \rightarrow \mathbb{T}^n$ defined by $x \mapsto x+\alpha$, and note that 
the $\widehat{r}$-orbit 
$\{n\alpha \co n\in \mathbb{Z}\}$ 
is a dense subset of a nontrivial subgroup of $\mathbb{T}^n$. 

 Recall that $h = [g^{-1}, s^{-1}]\circ[g^{-1}, s],$ with $s = r^M$ for $M$ as chosen above. 
 
 \begin{lemma} \label{small support lemma}
     With notation as defined above, the support of $h$ is contained in $X$.
 \end{lemma}   
 
 \begin{proof}
 Let $y \in [0,1) \setminus X$ be given; it will be shown that $h$ fixes $y$. Let $j, k \in \lbrace 1, \ldots, n \rbrace$ be such that $y \in I_j$ and $g(y) \in I_k$.
  By the definition of $X,$ $y$ is located a distance of at least $\epsilon$ away 
  from each of the endpoints $\beta_{j-1}$ and $\beta_j$ of $I_j$. Also, $y$ is at least $\epsilon$ away from any discontinuity of $g$; thus, the entire neighborhood $N_\epsilon(y)$
  is translated under $g$ by $\omega = g(y) - y$. Since $P$ also contains the points $g^{-1}(\beta_i)$ for $0 \leq i < n$, it follows that  
  $g(y)$ is located a distance of at least $\epsilon$ away from the endpoints of $I_k$. 
  
  With these conditions in mind, we trace the orbit of $y$ through the composition 
  \[ h = (g^{-1}s^{-1} g s)(g^{-1}sgs^{-1}).\]
Let $\gamma_j$ denote the rotation (mod $|I_j|$) induced by $s$ on $I_j$; define $\gamma_k$ similarly. By the construction of $s = r^M,$
we have $|\gamma_\iota| < \epsilon /10$ for $\iota \in \{j,k\}$. Thus, $s^{-1}(y) = y - \gamma_j$. As $s^{-1}(y)$ is still located in $N_\epsilon(y)$, 
we have $g(s^{-1}(y)) = y - \gamma_j + \omega = g(y) - \gamma_j$. Since $g(y) - \gamma_j$ is still at least a distance of $\epsilon / 2$ away from the endpoints
of $I_k$, we have $s(gs^{-1}(y)) = y - \gamma_j + \omega + \gamma_k = g(y) - \gamma_j + \gamma_k$. The map $g^{-1}$ translates $N_\epsilon(g(y))$ by $-\omega$, 
and $(g(y) - \gamma_j + \gamma_k) \in N_\epsilon(g(y))$. Thus $g^{-1}(sgs^{-1}(y)) =  y - \gamma_j + \gamma_k$.      
	
Let $z = g^{-1}sgs^{-1}(y)$, and note that $|y - z| < \epsilon / 5.$ Consequently, upon tracing the action of $g^{-1}s^{-1} g s$
on $z$, reasoning similar to the previous paragraph shows that $g^{-1}s^{-1} g s(z) = z + \gamma_j - \gamma_k = y.$ 
Thus, $h(y) = y$, as desired.
\end{proof}

If it happens that the map $h$ is the identity, then this suffices to show that $\langle r, g \rangle$ is not isomorphic to $F_2$. Suppose, however, 
that $h \neq id.$ 
Recall that $d$ is chosen so that $r^d(X')$ and $X'$ are disjoint.  Consider the interval exchange $k$ defined by 
\[ k = r^d h r^{-d}. \]
By Lemma \ref{small support lemma}, $h$ is supported in $X,$ and consequently, $k$ 
is supported in $r^d(X)$.  If the map $r$ has no periodic points, then $\text{supp}(r) = [0,1)$ and $X' = X$. Thus, in this situation 
the maps $h$ and $k$ have disjoint supports. It follows that $h$ and $k$ are commuting, nontrivial elements of $\langle r, g \rangle$ such that 
 $k \notin \langle h \rangle,$  which proves
that the group $\langle r, g \rangle$ is not isomorphic to $F_2$ when $r$ has no periodic points. 
	
To handle the general situation where $\text{Fix}(r) \neq \emptyset,$ it is shown below that the commutator $T = kh^{-1}k^{-1}h$ has
finite order that divides six. This again implies the existence of a nontrivial relation in $\langle r, g \rangle,$ completing the proof of Theorem 
\ref{disjoint rotation relation theorem}. The proof that $T^6 = id$ does not rely on the maps involved being IETs; it follows from the existence
of a conjugacy and the relation between the supports of $h$ and $k$. In the argument below, the \emph{support} of a map again refers to 
the complement of its set of fixed points. 

\begin{proposition}  \label{finite order commutator prop}
Let $h$ and $\phi$ be bijections of a set $\Omega.$  Write $\text{supp}(h) = A \sqcup B,$  where $A =  \text{supp}(h) \cap \text{supp}(\phi)$
and $B = \text{supp}(h) \cap \text{Fix}(\phi)$.  Let $k = \phi h \phi^{-1}$, and let $T = k h^{-1}k^{-1}h$. 
If $A \cap \phi(A) = \emptyset,$ then $T^6 = id$.
\end{proposition}

\begin{remark}
By Lemma \ref{small support lemma} and the condition that $X' \cap r^d(X') = \emptyset,$ it follows
 that the IETs $h$ and $\phi = r^d$ defined previously satisfy the hypotheses of Proposition \ref{finite order commutator prop}. 
\end{remark}
 
 \begin{proof}
 	We assume that $B$ is nonempty, since otherwise $h$ and $k$ have disjoint supports, in which case $T = id$.
For notation, let $C = \phi(A)$, and observe that $\text{supp}(k) = C \sqcup B.$ Since $k = \phi h \phi^{-1}$ and $\phi$ fixes 
all points in $B$, the following are true: 
\begin{itemize}
	\item[(a)] If $q \in B$ and $h^{\pm 1}(q) \in B$, then $k^{\pm 1}(q) = h^{\pm 1}(q)$.
	
	\item[(b)] If $q \in B$ and $h^{\pm 1}(q) \in A$, then $k^{\pm 1}(q) \in C$.
\end{itemize}
In addition, since $\text{supp}(h) = A \sqcup B$, one of (a) or (b) must hold for every $x \in B$.

  Note that if $p \in \Omega \setminus (A \sqcup B \sqcup C),$
then $T(p) = p.$ All other $p \in \Omega$ satisfy exactly one of the following situations:

\begin{itemize}
	\item[(I)] If $p \in A$ and $h(p) \in A$, then $T(p) = p$.  {\bf ($\ast$)}
	
	\item[(II)] If $p \in A$ and $h(p) \in B$, then $T(p) = h(p)$.  ${\bf (A \rightarrow B)}$
		
	\item[(III)] If $p \in C$ and $k^{-1}(p) \in C$, then $T(p) = p$. {\bf ($\ast$)}
	
	\item[(IV)] If $p \in C$ and $k^{-1}(p) \in B$, then one of the following holds:
			\begin{itemize}
				\item[(IVa)] $h^{-1}k^{-1}(p) \in A,$  in which case $T(p) = h^{-1}k^{-1}(p).$   ${\bf (C \rightarrow A)}$
				
				\item[(IVb)] $h^{-1}k^{-1}(p) \in B,$  in which case $T(p) = k^{-1}(p).$  ${\bf (C \rightarrow B)}$
			\end{itemize}
	
	\item[(V)] If $p \in B$ and $h(p) \in B$, then one of the following holds:
			\begin{itemize}
				\item[(Va)] $h^{-1}(p) \in A,$  in which case $T(p) = h^{-1}(p).$   ${\bf (B \rightarrow A)}$
				
				\item[(Vb)] $h^{-1}(p) \in B,$  in which case $T(p) = (p).$ {\bf ($\ast$)}
			\end{itemize}
	
	\item[(VI)] If $p \in B$ and $h(p) \in A$, then $T(p) = k(p).$   ${\bf (B \rightarrow C)}$		
\end{itemize}

The assertions in each of the above situations can be verified by tracing the location (either $A$, $B$, or $C$) of the
 point $p$ through the commutator $T = k h^{-1}k^{-1}h;$ in doing this, one uses the conditions (a) and (b) listed above, as 
 well as the fact that $h^{\pm 1} = id$ on $C$ and  $k^{\pm 1} = id$ on $A$.  
 
 For instance, situation (II) is checked as follows. 
 By assumption, $h(p) \in B$. Then, by (b) applied to $q = h(p)$, we have $k^{-1}h(p) \in C$. Since $h^{-1}$ acts trivially on $C$, 
 we have $h^{-1}k^{-1}h(p) = k^{-1}h(p)$, and thus $T(p) = k(h^{-1}k^{-1}h(p)) = h(p),$ as claimed. This trace can be summarized 
 in the following way: 
 \[  \left[ \begin{array}{c} p \\ A \\ \end{array} \right]         \begin{array}{c}  \mapsto \\ h \\ \end{array}    \left[ \begin{array}{c} h(p) \\  B \\ \end{array} \right]
 	 \begin{array}{c}  \mapsto \\ k^{-1} \\ \end{array}      \left[ \begin{array}{c}  k^{-1}h(p) \\ C \\ \end{array} \right]     
	 \begin{array}{c}  \mapsto  \\ h^{-1} \\ \end{array}      \left[ \begin{array}{c}  k^{-1}h(p) \\ C \\ \end{array} \right]     
 	\begin{array}{c}  \mapsto \\ k \\ \end{array}      \left[  \begin{array}{c}  h(p) \\ B \\ \end{array} \right]      \]  
 
 To verify situation (IVa) as a further example, the trace is illustrated by: 
  \[  \left[ \begin{array}{c} p \\ C \\ \end{array} \right]        \begin{array}{c}  \mapsto \\ h \\ \end{array}     \left[ \begin{array}{c} p \\  C \\ \end{array} \right]
 	 \begin{array}{c}  \mapsto \\ k^{-1} \\ \end{array}     \left[ \begin{array}{c}  k^{-1}(p) \\ B \\ \end{array} \right]    
	 \begin{array}{c}  \mapsto  \\ h^{-1} \\ \end{array}      \left[ \begin{array}{c}  h^{-1}k^{-1}(p) \\ A \\ \end{array} \right]    
 	\begin{array}{c}  \mapsto \\ k \\ \end{array}     \left[  \begin{array}{c}  h^{-1}k^{-1}(p) \\ A \\ \end{array} \right]      \]  
In situation (IVa), the application of $h$ is trivial since $h$ is the identity on $C$, the application of $k$ is trivial since $k$ is the identity on $A$, and all other 
labels and locations follow by the assumptions of (IVa). The remaining situations labelled by roman numerals can be verified through similar reasoning.

A few useful observations about $T$ are drawn from the above list of cases. First, for any $\Delta \in \lbrace A, B, C \rbrace$, if $p$ and $T(p)$ are 
both in $\Delta$, then $p$ is fixed by $T$; this corresponds to the situations (I), (III), and (Vb) that are marked by $(\ast)$ in the above list.
Further, note that for each ordered pair $(\Delta, \Gamma) \in \lbrace A, B, C \rbrace^2$ in which $\Delta \neq \Gamma$, there is at most one 
situation above such that $p \in \Delta$ and $T(p) \in \Gamma$; note that all such pairs are represented, with the exception of $(A \rightarrow C)$. 

These properties are now used to show that $T^6|_B = id_B.$ Let $x \in B$; if $T(x) \in B$, then $T(x) = x$. Otherwise, there are two cases to consider: \\

\noindent Case I: Suppose that $T(x) \in A$. Then $x$ must satisfy situation (Va), which implies that $h(x) \in B$ and $h^{-1}(x) \in A$. Thus $y = T(x) = h^{-1}(x).$ 
It is now the case that $y \in A$ and $h(y) = x \in B,$ so $y$ satisfies situation (II). It follows that $T(y) = h(y) = x,$ which shows that $T^2(x) = x$ in this case. \\

\noindent Case II: Suppose that $T(x) \in C$. Then $x$ satisfies situation (VI), so $h(x) \in A$ and $y = T(x) = k(x)$. Thus $y\in C$ and $k^{-1}(y) = x \in B$, so 
$y$ must satisfy either situation (IVa) or (IVb). If $y$ satisfies (IVb), then $T(y) = k^{-1}(y) = x$, in which case $T^2(x) = x$. Otherwise, assume that $y$
satisfies (IVa). Then $h^{-1}k^{-1}(y) = h^{-1}(x) \in A$ and $z = T(y) = h^{-1}k^{-1}(y) = h^{-1}(x) \in A$. Since $z \in A$ and $T(z) \neq z,$ it must be the case 
that $z$ satisfies situation (II). Thus, $T(z) = h(z) = x$, and it follows that $T^3(x) = x$. \\

As all points $x\in B$ have a $T$--orbit consisting of one, two, or three points, it follows that $T^6|_B = id_B.$ 
Considering $T^6|_A,$ note that points in $A \cap T^{-1}(A)$ are fixed by $T$, and note that points in $A \cap T^{-1}(B)$ are fixed by 
$T^6$ from the above argument. Since $T(A) \subset A\cup B$, it follows that $T^6|_A = id_A.$ Similarly, since points in $C \cap T^{-1}(C)$ 
are fixed by $T$ and all other points in $C$ map to $A \cup B$, we also have that $T^6|_C = id_C, $ which completes the proof that $T^6 = id.$ 
  \end{proof}

\bibliographystyle{abbrv}

\bibliography{Novak_ref_restrotns}

\begin{thebibliography}{10}

\bibitem{Arnoux81}
P.~Arnoux.
\newblock \'{E}changes d'intervalles et flots sur les surfaces.
\newblock In {\em Ergodic theory ({S}em., {L}es {P}lans-sur-{B}ex, 1980)
  ({F}rench)}, volume~29 of {\em Monograph. Enseign. Math.}, pages 5--38. Univ.
  Gen\`eve, Geneva, 1981.

\bibitem{AvilaForni07}
A.~Avila and G.~Forni.
\newblock Weak mixing for interval exchange transformations and translation
  flows.
\newblock {\em Ann. of Math. (2)}, 165(2):637--664, 2007.

\bibitem{AvilaViana07}
A.~Avila and M.~Viana.
\newblock Simplicity of {L}yapunov spectra: proof of the {Z}orich-{K}ontsevich
  conjecture.
\newblock {\em Acta Math.}, 198(1):1--56, 2007.

\bibitem{BrinSquier85}
M.~G. Brin and C.~C. Squier.
\newblock Groups of piecewise linear homeomorphisms of the real line.
\newblock {\em Invent. Math.}, 79(3):485--498, 1985.

\bibitem{delaHarpe00}
P.~de~la Harpe.
\newblock {\em Topics in geometric group theory}.
\newblock Chicago Lectures in Mathematics. University of Chicago Press,
  Chicago, IL, 2000.

\bibitem{FarbShalen02}
B.~Farb and P.~Shalen.
\newblock Groups of real-analytic diffeomorphisms of the circle.
\newblock {\em Ergodic Theory Dynam. Systems}, 22(3):835--844, 2002.

\bibitem{FerencziHoltonZamboni04}
S.~Ferenczi, C.~Holton, and L.~Q. Zamboni.
\newblock Structure of three-interval exchange transformations {III}: ergodic
  and spectral properties.
\newblock {\em J. Anal. Math.}, 93:103--138, 2004.

\bibitem{Ghys01}
{\'E}.~Ghys.
\newblock Groups acting on the circle.
\newblock {\em Enseign. Math. (2)}, 47(3-4):329--407, 2001.

\bibitem{Margulis00}
G.~Margulis.
\newblock Free subgroups of the homeomorphism group of the circle.
\newblock {\em C. R. Acad. Sci. Paris S\'er. I Math.}, 331(9):669--674, 2000.

\bibitem{Navas07bookEngArxiv}
A.~Navas.
\newblock Groups of circle diffeomorphisms ({E}nglish translation of {\em
  {G}rupos de difeomorfismos del c\'\i rculo}). 
\newblock {\tt arxiv:math/0607481v3 [math.DS] }

\bibitem{Navas07book}
A.~Navas.
\newblock {\em Grupos de difeomorfismos del c\'\i rculo}, volume~13 of {\em
  Ensaios Matem\'aticos [Mathematical Surveys]}.
\newblock Sociedade Brasileira de Matem\'atica, Rio de Janeiro, 2007.

\bibitem{Novak09}
C.~F. Novak.
\newblock Discontinuity-growth of interval-exchange maps.
\newblock {\em J. Mod. Dyn.}, 3(3):379--405, 2009.

\bibitem{Novak10}
C.~F. Novak.
\newblock Continuous interval exchange actions.
\newblock {\em to appear in Algebr. Geom. Topol.}, 2010.
\newblock {\tt arXiv:1007.1221v1 [math.DS] }

\bibitem{Tits72}
J.~Tits.
\newblock Free subgroups in linear groups.
\newblock {\em J. Algebra}, 20:250--270, 1972.

\bibitem{Veech84a}
W.~A. Veech.
\newblock The metric theory of interval exchange transformations. {I}.
  {G}eneric spectral properties.
\newblock {\em Amer. J. Math.}, 106(6):1331--1359, 1984.

\bibitem{Viana06}
M.~Viana.
\newblock Ergodic theory of interval exchange maps.
\newblock {\em Rev. Mat. Complut.}, 19(1):7--100, 2006.

\bibitem{Zorich06}
A.~Zorich.
\newblock Flat surfaces.
\newblock In {\em Frontiers in number theory, physics, and geometry. {I}},
  pages 437--583. Springer, Berlin, 2006.

\end{thebibliography}

\end{document}